\RequirePackage{fix-cm}
\RequirePackage{fixltx2e}
\documentclass[oneside,english]{amsart}

\usepackage[latin9]{inputenc}
\usepackage{babel}
\usepackage{url}
\usepackage{amsthm}
\usepackage{amssymb}
\usepackage{esint}
\usepackage[unicode=true,pdfusetitle,
 bookmarks=true,bookmarksnumbered=false,bookmarksopen=false,
 breaklinks=false,pdfborder={0 0 1},backref=false,colorlinks=false]
 {hyperref}
 \makeatletter

\renewcommand{\theequation}{\theequation. \arabic{equation}}
\numberwithin{equation}{section}
\newtheorem{thm}{Theorem}[section]

\newtheorem{prop}{Proposition}[section]

\usepackage[numbers,sort&compress]{natbib}

\def\squarebox#1{\hbox to #1{\hfill\vbox to #1{\vfill}}}
\def\qed{\hspace*{\fill}
         \vbox{\hrule\hbox{\vrule\squarebox{.667em}\vrule}\hrule}\smallskip}
\begin{document}\large
\title[Appell series $F_{1}$ over finite fields]
{\Large Another finite field analogue for Appell series $F_{1}$}
\author{\small  Bing He}
\address{\small
College of Science, Northwest A\&F University,
   Yangling 712100, Shaanxi, People's Republic of China}
\email{yuhe001@foxmail.com; yuhelingyun@foxmail.com}




\begin{abstract}
\small In this paper we introduce another finite field analogue for Appell series $F_{1}$ and
obtain certain  reduction formulae and  a generating function for this analogue.
\end{abstract}
\maketitle
\section{Introduction}
Let $q$ be a power of a prime and let $\mathbb{F}_{q}$ and $\widehat{\mathbb{F}}^{*}_{q}$  denote the finite field of $q$ elements and  the group of multiplicative characters of $\mathbb{F}^{*}_{q}$ respectively. Then the domain of all characters $\chi$ of $\mathbb{F}^{*}_{q}$ can be extended to $\mathbb{F}_{q}$ by setting $\chi(0)=0$ for all characters. $\overline{\chi}$ and $\varepsilon$ are denoted as  the inverse of $\chi$ and the trivial character respectively.  See \cite{BEW} and \cite [Chapter 8]{IR} for more information about characters.

Define the generalized hypergeometric function by \cite{B}
$$ {}_{n+1}F_n \left(\begin{matrix}
a_0, a_1, \ldots , a_{n} \\
b_1, \ldots , b_n \end{matrix}
\bigg| x \right):=\sum_{k=0}^{\infty}\frac{(a_{0})_{k}(a_{1})_{k}\cdots(a_{n})_{k}}{k!(b_{1})_{k}\cdots(b_{n})_{k}}x^{k}$$
where $(z)_{k}$ is the Pochhammer symbol given by
\begin{equation*}
(z)_{0}=1,~(z)_{k}=z(z+1)\cdots(z+k-1)\text{ for } k\geq 1.
\end{equation*}

It was Greene  who in \cite{Gr}  developed the theory of hypergeometric functions over finite fields and derived various transformation and summation identities for hypergeometric series over finite fields which are finite field analogues of those in the classical case.  In that paper, the notation
\begin{equation*}
{}_{2}F_1 \left(\begin{matrix}
A, B \\
C \end{matrix}
\bigg| x \right)^{G}=\varepsilon(x)\frac{BC(-1)}{q}\sum_{y}B(y)\overline{B}C(1-y)\overline{A}(1-xy)
\end{equation*}
for $A,B,C\in \widehat{\mathbb{F}}_{q}$ and $x\in \mathbb{F}_{q},$ which is
the finite field analogue of the integral representation of Gauss hypergeometric series \cite{B}:
\begin{equation*}
{}_{2}F_1 \left(\begin{matrix}
a, b \\
c \end{matrix}
\bigg| x \right)=\frac{\Gamma(c)}{\Gamma(b)\Gamma(c-b)}\int_{0}^{1}t^b(1-t)^{c-b}(1-tx)^{-a}\frac{dt}{t(1-t)},
\end{equation*}
was introduced by Greene and the finite field analogue of the binomial coefficient:
\begin{equation*}
  {A\choose B}^{G}=\frac{B(-1)}{q}J(A,\overline{B}),
\end{equation*}
where $J(\chi,\lambda)$ is the Jacobi sum given by $$J(\chi,\lambda)=\sum_{u}\chi(u)\lambda(1-u),$$ was also defined.  For more information about the finite field analogue of the generalized hypergeometric functions, please see \cite{FL, M, EG}.

In this paper, for the sake of simplicity, we use the notation
\begin{equation*}
  {A\choose B}=q{A\choose B}^{G}=B(-1)J(A,\overline{B}).
\end{equation*}
Furthermore, we define the finite field analogue of the classic Gauss hypergeometric series as
\begin{equation*}
{}_{2}F_1 \left(\begin{matrix}
A, B \\
C \end{matrix}
\bigg| x \right)=q\cdot{}_{2}F_1 \left(\begin{matrix}
A, B \\
C \end{matrix}
\bigg| x \right)^{G}=\varepsilon(x)BC(-1)\sum_{y}B(y)\overline{B}C(1-y)\overline{A}(1-xy).
\end{equation*}
Then by \cite [Theorem 3.6]{Gr},
\begin{equation*}
  {}_{2}F_1 \left(\begin{matrix}
A, B \\
C \end{matrix}
\bigg| x \right)=\frac{1}{q-1}\sum_{\chi}{A\chi\choose \chi}{B\chi\choose C\chi}\chi(x)
\end{equation*}
for any $A,B,C\in \widehat{\mathbb{F}}_{q}$ and $x\in \mathbb{F}_{q}.$  Similarly, the finite field analogue of the generalized hypergeometric series for any $A_{0},A_{1},\cdots, A_{n},B_{1},\cdots,B_{n}\in \widehat{\mathbb{F}}_{q}$ and $x\in \mathbb{F}_{q}$ is defined by
\begin{equation*}
  {}_{n+1}F_n \left(\begin{matrix}
A_{0}, A_{1},\cdots, A_{n} \\
 B_{1},\cdots, B_{n} \end{matrix}
\bigg| x \right)=\frac{1}{q-1}\sum_{\chi}{A_{0}\chi\choose \chi}{A_{1}\chi\choose B_{1}\chi}\cdots{A_{n}\chi\choose B_{n}\chi}\chi(x).
\end{equation*}
In our notations one of Greene's theorems is as follows.
\begin{thm}\emph{(See \cite [Theorem 4.9]{Gr})}For any characters $A,B,C\in \widehat{\mathbb{F}}_{q},$  we have
\begin{equation}\label{e1-1}
{}_{2}F_1 \left(\begin{matrix}
A, B \\
 C \end{matrix}
\bigg| 1 \right)=A(-1){B\choose \overline{A}C}.
\end{equation}
\end{thm}

Among these interesting double hypergeometric functions in the field of hypergeometric functions, Appell's four functions may be the most important functions. See \cite{A, B, CA, S} for more detailed material about Appell's functions. One of them is as follows:
\begin{align*}
F_{1}(a;b,b';c;x,y)&=\sum_{m,n\geq 0}\frac{(a)_{m+n}(b)_{m}(b')_{n}}{m!n!(c)_{m+n}}x^my^n,~|x|<1,~|y|<1.
\end{align*}

Inspired by Greene's work, Li \emph{et al} in \cite{LLM} gave a finite field analogue of the Appell series $F_{1}$ and
obtained some transformation and reduction formulas and the generating functions for the function over finite fields.  In that paper, the finite field analogue of the Appell series $F_{1}$ was given by
\begin{equation*}
F_{1}(A;B,B';C;x,y)=\varepsilon(xy)AC(-1)\sum_{u}A(u)\overline{A}C(1-u)\overline{B}(1-ux)\overline{B'}(1-uy)
\end{equation*}
owing to the fact that the $F_{1}$ function has the integral representation in terms of
a single integral \cite [Chapter IX]{B}:
\begin{equation*}
F_{1}(a;b,b';c;x,y)=\frac{\Gamma(c)}{\Gamma(a)\Gamma(c-a)}\int_{0}^{1}u^{a-1}(1-u)^{c-a-1}(1-ux)^{-b}(1-uy)^{-b'}du
\end{equation*}
where $0< \Re(a)<\Re(c).$

Motivated by the work of Greene \cite{Gr} and the second author \emph{et al} \cite{LLM}, we give another finite field analogue of  the Appell series  $F_{1}.$  Since the Appell series  $F_{1}$ has the following double integral representation  \cite [Chapter IX]{B}:
\begin{align*}
F_{1}(a;b,b';c;x,y)&=\frac{\Gamma(c)}{\Gamma(b)\Gamma(b')\Gamma(c-b-b')}\\
&\cdot \int\int u^{b-1}v^{b'-1}(1-u-v)^{c-b-b'-1}(1-ux-vy)^{-a}dudv
\end{align*}
where the sums are taken over the triangle region $\{(u,v)|u\geq 0, v\geq 0,u+v\leq 1\},$
we now give the  finite field analogue of $F_{1}$ in the following form:
\begin{align*}
  F_{1}(A;B,B';C;x,y)=\varepsilon(xy)BB'(-1)\sum_{u,v}B(u)B'(v)C\overline{B}\overline{B'}(1-u-v)\overline{A}(1-ux-vy),
\end{align*}
where $A,B,B',C \in \widehat{\mathbb{F}}_{q},~x,y\in \mathbb{F}_{q}$ and  each sum ranges over all the elements of $\mathbb{F}_{q}.$ In the above definition,  the factor $\frac{\Gamma(c)\Gamma(c')}{\Gamma(b)\Gamma(b')\Gamma(c-b)\Gamma(c'-b')}$ is dropped to obtain simpler results. We choose the factor  $\varepsilon(xy)\cdot BB'CC'(-1)$ to get a better expression in terms of binomial coefficients.

From the definition of $F_{1}(A;B,B';C;x,y)$ we know that
for any $A,B,B',C,C'\in \widehat{\mathbb{F}}_{q}$ and $x,y\in \mathbb{F}_{q},$ we have
\begin{align}
  F_{1}(A;B,B';C;x,y)&=F_{1}(A;B',B;C;y,x). \label{e1-2}
\end{align}

The aim of this paper is to give certain reduction formulas and  a generating function for the Appell series $F_{1}$ over finite fields.
The fact that the Appell series $F_{1}$  has a  double one leads us to giving  a finite field analogue for the Appell series $F_{1}$ which is more complicated than that for $F_{1}$ in \cite{LLM}. Consequently,  the results on  reduction formulas and that generating function for this analogue are also  more complicated than those in \cite{LLM}.

Another expression for $F_{1}(A;B,B';C;x,y)$ will be given in the next section. We give certain reduction formulae for $F_{1}(A;B,B';C;x,y)$ in Section 3. The last section is devoted to deducing a generating function for $F_{1}(A;B,B';C;x,y).$
\section{Another expression for $F_{1}(A;B,B';C;x,y)$}
In this section we present another expression for $F_{1}(A;B,B';C;x,y).$
\begin{thm}\label{t2-1} For $A,B,B',C \in \widehat{\mathbb{F}}_{q},~x,y\in \mathbb{F}_{q},$ we have
\begin{align*}
F_{1}(A;B,B';C;x,y)&=\frac{1}{(q-1)^2}\sum_{\lambda,\mu}{A\lambda\choose \lambda}{A\lambda\mu\choose \mu}{\overline{B}\overline{B'}C \choose \overline{B'}C\lambda}{\overline{B'}C\lambda\choose C\lambda\mu}\lambda(-x)\mu(-y)\notag\\
&~+C(-1){AB'\overline{C}\choose B'\overline{C}}B'\overline{C}(x)\overline{A}\overline{B'}C(1-y)\notag\\
&~+BC(-1)\overline{A}(x)A\overline{C}(y){AB'\overline{C}\choose A\overline{B}}\overline{A}\overline{B'}C(1-y)\\
&~+(q-1)\overline{A}(x)B'(-1)\delta(A\overline{C}B')\delta(y-1),
\end{align*}
where each sum ranges over all multiplicative characters of $\mathbb{F}_{q}.$
\end{thm}
To carry out our study, we need some auxiliary results which will be used  in the sequel.

\begin{prop}\label{pp2-1} \emph{(See \cite [(2.6), (2.7), (2.8) and (2.12)]{Gr})} If $A,B\in \widehat{\mathbb{F}_{q}},$ then
\begin{align}
  {A\choose B}&={A\choose A\overline{B}},\label{f2}\\
  {A\choose B}&={B\overline{A}\choose B}B(-1)\label{f1},\\
  {A\choose B}&={\overline{B}\choose \overline{A}}AB(-1),\label{f3}\\
{A\choose \varepsilon}&={A\choose A}=-1+(q-1)\delta(A),\label{f4}
\end{align}
where $\delta(\chi)$ is given by
\begin{equation*}
\delta(\chi)=\left\{
               \begin{array}{ll}
                 1  & \hbox{if $\chi=\varepsilon$} \\
                 0  & \hbox{otherwise}
               \end{array}.
             \right.
\end{equation*}
\end{prop}

\begin{prop}\emph{(Binomial theorem, see \cite [(2.5)]{Gr})} For any character $A\in \widehat{\mathbb{F}}_{q}$ and $x\in \mathbb{F}_{q},$ we have
\begin{equation*}
  A(1+x)=\delta(x)+\frac{1}{q-1}\sum_{\chi}{A\choose \chi}\chi(x),
\end{equation*}
where the sum ranges over all multiplicative characters of $\mathbb{F}_{q}$ and $\delta(x)$ is a function on $\mathbb{F}_{q}$ given by
\begin{equation*}
\delta(x)=\left\{
            \begin{array}{ll}
              1 & \hbox{if $x=0$} \\
              0 & \hbox{if $x\neq 0$}
            \end{array}.
          \right.
\end{equation*}
\end{prop}

We are now ready to show Theorem \ref{t2-1}.\\
\emph{Proof of Theorem \ref{t2-1}.} It is clear that $F_{1}(A;B,B';C;x,y)=0$ for $y=0.$ We now consider the case $y\neq 0.$
 When $v\neq 1,$ it is known from the binomial theorem that
\begin{equation*}
  C\overline{B}\overline{B'}\left(1-\frac{u}{1-v}\right)=\delta(u)+\frac{1}{q-1}\sum_{\chi}{C\overline{B}\overline{B'}\choose \chi}\chi(-u)\overline{\chi}(1-v).
\end{equation*}
Then
\begin{align}
\label{t2-1-1} C\overline{B}\overline{B'}(1-u-v)&=C\overline{B}\overline{B'}(1-v)C\overline{B}\overline{B'}\left(1-\frac{u}{1-v}\right)\\
 &=\delta(u)C\overline{B}\overline{B'}(1-v)+\frac{1}{q-1}\sum_{\chi}{C\overline{B}\overline{B'}\choose \chi}\chi(-u)C\overline{B}\overline{B'}\overline{\chi}(1-v).  \notag
\end{align}
Similarly, when $v\neq y^{-1},$
\begin{align}\label{t2-1-2}
\overline{A}(1-ux-vy)=\delta(ux)\overline{A}(1-vy)+\frac{1}{q-1}\sum_{\lambda}{\overline{A}\choose \lambda}\lambda(-ux)\overline{A}\overline{\lambda}(1-vy).
\end{align}
In particular, when $y=1$ and $v\neq 1,$
\begin{equation}\label{t2-1-3}
\overline{A}(1-ux-v)=\delta(ux)\overline{A}(1-v)+\frac{1}{q-1}\sum_{\lambda}{\overline{A}\choose \lambda}\lambda(-ux)\overline{A}\overline{\lambda}(1-v).
\end{equation}
If $y=1,$ then substituting  \eqref{t2-1-1}, \eqref{t2-1-3} and combining the fact that $\varepsilon(x)B(u)\delta(ux)=B(u)\delta(u)=0,$ \cite [(1.15)]{Gr}, \eqref{f2}--\eqref{f3}  gives
\begin{align*}
&\varepsilon(x)\sum_{\scriptstyle u\in \mathbb{F}_{q} \atop \scriptstyle v\neq 1}B(u)B'(v)C\overline{B}\overline{B'}(1-u-v)\overline{A}(1-ux-v)\\
&=\frac{1}{(q-1)^2}\sum_{\chi,\lambda}\chi(-1)\lambda(-x){C\overline{B}\overline{B'}\choose \chi}{\overline{A}\choose \lambda}\sum_{v}B'(v)C\overline{A}\overline{B}\overline{B'}\overline{\chi}\overline{\lambda}(1-v)\sum_{u}B\chi\lambda(u)\\
&=\frac{B(-1)}{q-1}\sum_{\lambda}{C\overline{B}\overline{B'}\choose \overline{B}\overline{\lambda}}{A\lambda\choose \lambda}\lambda(-x)\sum_{v}B'(v)C\overline{A}\overline{B'}(1-v)\\
&=\frac{AB(-1)}{q-1}{B'\choose AB'\overline{C}}\sum_{\lambda}{A\lambda\choose \lambda}{B\lambda\choose C\overline{B'}\lambda}\lambda(x)\\
&=AB(-1){B'\choose \overline{A}C}{}_{2}F_1 \left(\begin{matrix}
A, B \\
C\overline{B'} \end{matrix}
\bigg| x \right).
\end{align*}
On the other hand, by \cite [(1.15)]{Gr},
\begin{align*}
\varepsilon(x)\sum_{u\in \mathbb{F}_{q}}B(u)C\overline{B}\overline{B'}(-u)\overline{A}(-ux)&=ABB'C(-1)\overline{A}(x)\sum_{u\in \mathbb{F}_{q}}\overline{A}\overline{B'}C(u)\\
&=(q-1)\overline{A}(x)B(-1)\delta(A\overline{C}B').
\end{align*}
Then
\begin{align}
 \label{t21-1}F_{1}(A;B,B';C;x,1)&=\varepsilon(x)BB'(-1)\left(\sum_{\scriptstyle u\in \mathbb{F}_{q} \atop \scriptstyle v\neq 1}+\sum_{\scriptstyle u\in \mathbb{F}_{q} \atop \scriptstyle v=1}\right)\\
 &=AB'(-1){B'\choose \overline{A}C}{}_{2}F_1 \left(\begin{matrix}
A, B \\
C\overline{B'} \end{matrix}
\bigg| x \right)\notag\\
&~+(q-1)\overline{A}(x)B'(-1)\delta(A\overline{C}B').\notag
\end{align}
If $y\neq 1,$ then we substitute \eqref{t2-1-2} and combine $\varepsilon(x)\delta(ux)B(u)=0,$ \eqref{f1} and \cite [(1.15)]{Gr} to get
\begin{align*}
&\varepsilon(x)\sum_{u\in \mathbb{F}_{q}}B(u)C\overline{B}\overline{B'}(-u)\overline{A}(1-ux-y)\\
&=\frac{1}{q-1}BB'C(-1)\sum_{\lambda}{A\lambda\choose \lambda}\lambda(x)\overline{A}\overline{\lambda}(1-y)\sum_{u\in \mathbb{F}_{q}}C\overline{B'}\lambda(u)\\
&=B(-1){AB'\overline{C}\choose B'\overline{C}}B'\overline{C}(-x)\overline{A}C\overline{B'}(1-y).
\end{align*}
Similarly, we combine \eqref{t2-1-1}, $B(u)\delta(u)=0,$  \eqref{f1} and \cite [(1.15)]{Gr} to obtain
\begin{align*}
&\varepsilon(x)\sum_{u}B(u)B'(y^{-1})C\overline{B}\overline{B'}(1-u-y^{-1})\overline{A}(-ux)\\
&=\frac{\overline{A}(-x)\overline{B'}(y)}{q-1}\sum_{\chi}{BB'\overline{C}\chi\choose \chi}C\overline{B}\overline{B'}\overline{\chi}(1-y^{-1})\sum_{u}\overline{A}B\chi(u)\\
&=\overline{A}(-x)A\overline{C}(y){AB'\overline{C}\choose A\overline{B}}\overline{A}\overline{B'}C(y-1).
\end{align*}
Again we combine \eqref{t2-1-1}, \eqref{t2-1-2}, $\varepsilon(x)B(u)\delta(ux)=B(u)\delta(u)=B'(v)\delta(v)=0,$  \eqref{f1}, the binomial theorem and \cite [(1.15)]{Gr} to derive
\begin{align*}
&\varepsilon(x)\sum_{\scriptstyle u\in \mathbb{F}_{q} \atop \scriptstyle v\neq 1~\text{and}~y^{-1}}B(u)B'(v)C\overline{B}\overline{B'}(1-u-v)\overline{A}(1-ux-vy)\\
&=\frac{1}{(q-1)^2}\sum_{\chi,\lambda}{BB'\overline{C}\chi\choose \chi}{A\lambda\choose \lambda}\lambda(x)\sum_{\scriptstyle u\in \mathbb{F}_{q} \atop \scriptstyle v\neq 1~\text{and}~y^{-1}}B\chi\lambda(u)B'(v)C\overline{B}\overline{B'}\overline{\chi}(1-v)\overline{A}\overline{\lambda}(1-vy)\\
&=\frac{1}{(q-1)^2}\sum_{\chi,\lambda}{BB'\overline{C}\chi\choose \chi}{A\lambda\choose \lambda}\lambda(x)\sum_{u,v\in \mathbb{F}_{q}}B\chi\lambda(u)B'(v)C\overline{B}\overline{B'}\overline{\chi}(1-v)\overline{A}\overline{\lambda}(1-vy)
\end{align*}
\begin{align*}
&=\frac{1}{(q-1)^3}\sum_{\chi,\lambda,\mu}{BB'\overline{C}\chi\choose \chi}{A\lambda\choose \lambda}{A\lambda\mu\choose \mu}\lambda(x)\mu(y)\sum_{v}B'\mu(v)C\overline{B}\overline{B'}\overline{\chi}(1-v)\sum_{u}B\chi\lambda(u)\\
&=\frac{1}{(q-1)^2}\sum_{\lambda,\mu}{B'\overline{C}\overline{\lambda}\choose \overline{B}\overline{\lambda}}{A\lambda\choose \lambda}{A\lambda\mu\choose \mu}\lambda(x)\mu(y)\sum_{v}B'\mu(v)C\overline{B'}\lambda(1-v)\\
&=\frac{BB'(-1)}{(q-1)^2}\sum_{\lambda,\mu}{\overline{B}\overline{B'}C \choose \overline{B'}C\lambda}{A\lambda\choose \lambda}{A\lambda\mu\choose \mu}{\overline{B'}C\lambda\choose C\lambda\mu}\lambda(-x)\mu(-y).
\end{align*}
Then
\begin{align}
\label{t21-2}F_{1}(A;B,B';C;x,y)&=\varepsilon(x)BB'(-1)\left(\sum_{\scriptstyle u\in \mathbb{F}_{q} \atop \scriptstyle v=1}+\sum_{\scriptstyle u\in \mathbb{F}_{q} \atop \scriptstyle v=y^{-1}}+\sum_{\scriptstyle u\in \mathbb{F}_{q} \atop \scriptstyle v\neq 1~\text{and}~y^{-1}}\right)\\
&=\frac{1}{(q-1)^2}\sum_{\lambda,\mu}{\overline{B}\overline{B'}C \choose \overline{B'}C\lambda}{A\lambda\choose \lambda}{A\lambda\mu\choose \mu}{\overline{B'}C\lambda\choose C\lambda\mu}\lambda(-x)\mu(-y)\notag\\
&~+C(-1){AB'\overline{C}\choose B'\overline{C}}B'\overline{C}(x)\overline{A}\overline{B'}C(1-y)\notag\\
&~+BC(-1)\overline{A}(x)A\overline{C}(y){AB'\overline{C}\choose A\overline{B}}\overline{A}\overline{B'}C(1-y). \notag
\end{align}
In addition, by \eqref{e1-1}
\begin{align}
&\frac{1}{(q-1)^2}\sum_{\lambda,\mu}{\overline{B}\overline{B'}C \choose \overline{B'}C\lambda}{A\lambda\choose \lambda}{A\lambda\mu\choose \mu}{\overline{B'}C\lambda\choose C\lambda\mu}\lambda(-x)\mu(-1)\label{t21-3}\\
&=\frac{B'(-1)}{(q-1)^2}\sum_{\lambda}{B\lambda \choose \overline{B'}C\lambda}{A\lambda\choose \lambda}\lambda(-x)\sum_{\mu}{A\lambda\mu\choose \mu}{B'\mu\choose C\lambda\mu}\mu(1)\notag\\
&=\frac{B'(-1)}{q-1}\sum_{\lambda}{B\lambda \choose \overline{B'}C\lambda}{A\lambda\choose \lambda}\lambda(-x){}_{2}F_1 \left(\begin{matrix}
A\lambda, B' \\
C\lambda \end{matrix}
\bigg| 1 \right)\notag\\
&=\frac{AB'(-1)}{q-1}{B'\choose \overline{A}C}\sum_{\lambda}{A\lambda\choose \lambda}{B\lambda \choose \overline{B'}C\lambda}\lambda(x)\notag\\
&=AB'(-1){B'\choose \overline{A}C}{}_{2}F_1 \left(\begin{matrix}
A, B \\
C\overline{B'} \end{matrix}
\bigg| x \right). \notag
\end{align}
In view of \eqref{t21-1}--\eqref{t21-3}, we get
\begin{align*}
F_{1}(A;B,B';C;x,y)&=\frac{1}{(q-1)^2}\sum_{\lambda,\mu}{\overline{B}\overline{B'}C \choose \overline{B'}C\lambda}{A\lambda\choose \lambda}{A\lambda\mu\choose \mu}{\overline{B'}C\lambda\choose C\lambda\mu}\lambda(-x)\mu(-y)\\
&~+C(-1){AB'\overline{C}\choose B'\overline{C}}B'\overline{C}(x)\overline{A}\overline{B'}C(1-y)
\end{align*}
\begin{align*}
&~+BC(-1)\overline{A}(x)A\overline{C}(y){AB'\overline{C}\choose A\overline{B}}\overline{A}\overline{B'}C(1-y)\\
&~+(q-1)\overline{A}(x)B'(-1)\delta(A\overline{C}B')\delta(y-1),
\end{align*}
which completes the proof of Theorem \ref{t2-1}.\qed
\section{Reduction formulae}
In this section we give certain reduction  formulae for $F_{1}(A;B,B';C;x,y).$
In order to derive these formulae we need some auxiliary results.
\begin{prop}\label{p2}\emph{(See \cite [Corollary 3.16 and Theorem 3.15]{Gr})} For any $A,B,C,D\in \widehat{\mathbb{F}_{q}}$ and $x \in \mathbb{F}_{q},$ we have
\begin{align}
  {}_{2}F_1 \left(\begin{matrix}
A, \varepsilon \\
C \end{matrix}
\bigg| x \right)&={C\choose A}A(-1)\overline{C}(x)\overline{A}C(1-x)
-C(-1)\varepsilon(x)\label{p2-1}\\
&~~+(q-1)A(-1)\delta(1-x)\delta(\overline{A}C),\notag\\
{}_{2}F_1 \left(\begin{matrix}
A, B \\
A \end{matrix}
\bigg| x \right)&={B\choose A}\varepsilon(x)\overline{B}(1-x)-\overline{A}(-x)\label{p2-2}\\
&~~+(q-1)A(-1)\delta(1-x)\delta(B),\notag\\
{}_{3}F_2 \left(\begin{matrix}
A,B,C\\
A,D \end{matrix}
\bigg| x \right)&={B\choose A}{}_{2}F_1 \left(\begin{matrix}
B,C \\
D \end{matrix}
\bigg| x \right)-\overline{A}(-x){C\overline{A}\choose D\overline{A}}\label{p2-3}\\
&~~+(q-1)A(-1)\overline{D}(x)\overline{C}D(1-x)\delta(B).\notag
\end{align}
\end{prop}

From the definition of $F_{1}(a;b,b';c;x,y)$ we know that
\begin{align}
  F_{1}(a;b,0;c;x,y)&={}_{2}F_1 \left(\begin{matrix}
a, b \\
c \end{matrix}
\bigg| x \right),\label{e3-1}\\
F_{1}(a;0,b';c;x,y)&={}_{2}F_1 \left(\begin{matrix}
a, b' \\
c \end{matrix}
\bigg| y \right). \label{e3-2}
\end{align}
We now give a finite field analogue for \eqref{e3-1}.
\begin{thm}\label{t31}Let $A,B, C\in \widehat{\mathbb{F}}_{q}$ and $x\in \mathbb{F}_{q},~y \in \mathbb{F}^{*}_{q}\backslash\{1\}.$  Then
\begin{align*}
  &F_{1}(A;B,\varepsilon;C;x,y)\\
&=-C(-1){}_{2}F_1 \left(\begin{matrix}
A, B \\
C \end{matrix}
\bigg| x \right)+\varepsilon(x)B\overline{C}(y)\overline{A}C(1-y)\overline{B}(y-x){A\overline{C}\choose A}{B\choose C}\\
&~+(q-1)\overline{C}(y)\overline{A}C(y-1)\delta(y-x)\delta(B){C\choose A}\\
&~+(q-1)\overline{C}(x)B\overline{C}(y)\overline{A}C(1-y)\overline{B}C(y-x)\delta(C).
\end{align*}
\end{thm}
\noindent{\it Proof.} It is easily seen from \eqref{p2-2} that
\begin{align*}
{}_{2}F_1 \left(\begin{matrix}
C,B\\
C \end{matrix}
\bigg| \frac{x}{y} \right)={B\choose C}\varepsilon(x)\overline{B}(y-x)B(y)-\overline{C}(-x)C(y)+(q-1)C(-1)\delta(y-x)\delta(B),
\end{align*}
which combines \eqref{p2-3} to  deduce
\begin{align}
 \label{t31-1}\frac{1}{q-1}\sum_{\lambda}{A\lambda\choose \lambda}{C\lambda\choose A\lambda}{B\lambda \choose C\lambda}\lambda\left(\frac{x}{y}\right)&={}_{3}F_2 \left(\begin{matrix}
A,C,B\\
A,C \end{matrix}
\bigg| \frac{x}{y} \right)\\
&={C\choose A}{}_{2}F_1 \left(\begin{matrix}
C,B\\
C \end{matrix}
\bigg| \frac{x}{y} \right)-\overline{A}(-x)A(y){B\overline{A}\choose C\overline{A}}\notag\\
&~+(q-1)A(-1)\overline{C}(x)B(y)\overline{B}C(y-x)\delta(C)\notag\\
&={C\choose A}{B\choose C}\varepsilon(x)\overline{B}(y-x)B(y)-{C\choose A}\overline{C}(-x)C(y)\notag\\
&~-\overline{A}(-x)A(y){B\overline{A}\choose C\overline{A}}+(q-1)C(-1){C\choose A}\delta(y-x)\delta(B)\notag\\
&~+(q-1)A(-1)\overline{C}(x)B(y)\overline{B}C(y-x)\delta(C).\notag
\end{align}
It follows from  \eqref{p2-1} that
\begin{align*}
\frac{1}{q-1}\sum_{\mu}{A\lambda\mu\choose \mu}{\mu\choose C\lambda\mu}\mu(y)&={}_{2}F_1 \left(\begin{matrix}
A\lambda, \varepsilon \\
C\lambda \end{matrix}
\bigg| y \right)\\
&={C\lambda\choose A\lambda}A\lambda(-1)\overline{C\lambda}(y)\overline{A}C(1-y)
-C\lambda(-1).
\end{align*}
Then
\begin{align}
&\frac{1}{(q-1)^2}\sum_{\lambda}{A\lambda\choose \lambda}{B\lambda \choose C\lambda}\lambda(-x)\sum_{\mu}{A\lambda\mu\choose \mu}{\mu\choose C\lambda\mu}\mu(y)\label{t31-2}\\
&=\frac{A(-1)\overline{C}(y)\overline{A}C(1-y)}{q-1}\sum_{\lambda}{A\lambda\choose \lambda}{C\lambda\choose A\lambda}{B\lambda \choose C\lambda}\lambda\left(\frac{x}{y}\right)\notag\\
&~-\frac{C(-1)}{q-1}\sum_{\lambda}{A\lambda\choose \lambda}{B\lambda \choose C\lambda}\lambda(x). \notag
\end{align}
Therefore, by \eqref{f1}, \eqref{t31-1}, \eqref{t31-2} and cancelling some terms, we have
\begin{align*}
F_{1}(A;B,\varepsilon;C;x,y)&=\frac{1}{(q-1)^2}\sum_{\lambda}{A\lambda\choose \lambda}{B\lambda \choose C\lambda}\lambda(-x)\sum_{\mu}{A\lambda\mu\choose \mu}{\mu\choose C\lambda\mu}\mu(y)\notag\\
&~+C(-1){A\overline{C}\choose \overline{C}}\overline{C}(x)\overline{A}C(1-y)+BC(-1)\overline{A}(x)A\overline{C}(y){A\overline{C}\choose A\overline{B}}\overline{A}C(1-y)\\
&=-C(-1){}_{2}F_1 \left(\begin{matrix}
A, B \\
C \end{matrix}
\bigg| x \right)+\varepsilon(x)B\overline{C}(y)\overline{A}C(1-y)\overline{B}(y-x){A\overline{C}\choose A}{B\choose C}\\
&~+(q-1)\overline{C}(y)\overline{A}C(y-1)\delta(y-x)\delta(B){C\choose A}\\
&~+(q-1)\overline{C}(x)B\overline{C}(y)\overline{A}C(1-y)\overline{B}C(y-x)\delta(C).
\end{align*}
This concludes the proof of Theorem \ref{t31}.\qed

From Theorem \ref{t31} and \eqref{e1-2} we can easily deduce a finite field analogue for \eqref{e3-2}.
\begin{thm}\label{t32}Let $A,B', C\in \widehat{\mathbb{F}}_{q}$ and $x\in \mathbb{F}^{*}_{q}\backslash\{1\},~y\in \mathbb{F}_{q}.$  Then
\begin{align*}
  &F_{1}(A;\varepsilon,B';C;x,y)\\
&=-C(-1){}_{2}F_1 \left(\begin{matrix}
A, B' \\
C \end{matrix}
\bigg| y \right)+\varepsilon(y)B'\overline{C}(x)\overline{A}C(1-x)\overline{B'}(x-y){A\overline{C}\choose A}{B'\choose C}\\
&~+(q-1)\overline{C}(x)\overline{A}C(x-1)\delta(x-y)\delta(B'){C\choose A}\\
&~+(q-1)\overline{C}(y)B'\overline{C}(x)\overline{A}C(1-x)\overline{B'}C(x-y)\delta(C).
\end{align*}
\end{thm}


\section{A generating function}
In this section, we establish a generating function for $F_{2}(A;B,B';C,C';x,y).$
\begin{thm}\label{t4-1}For any $A,B,B',C,C'\in \widehat{\mathbb{F}_{q}} $ and $x\in \mathbb{F}^*_{q},~y,t\in \mathbb{F}^*_{q}\backslash\{1\},$ we have
\begin{align*}
  &\frac{1}{q-1}\sum_{\theta}{A\theta\choose\theta}F_{1}(A\theta;B,B';C;x,y)\theta(t)\\
&=\overline{A}(1-t)F_{1}\left(A;B,B';C;\frac{x}{1-t}, \frac{y}{1-t}\right)\notag\\
&~-BC(-1)B'(1-t)\overline{A}(x)A\overline{C}(y){AB'\overline{C}\choose A\overline{B}}\overline{A}\overline{B'}C(1-t-y)\\
&~-\overline{A}(-t)B'(-1)\overline{C}(y)\overline{B'}C(1-y){}_{2}F_1 \left(\begin{matrix}
A, B \\
\overline{B'}C \end{matrix}
\bigg| \frac{x(1-y)}{ty}\right)\\
&~-\overline{A}(-t){\overline{B'}C\choose C}{\overline{B}\overline{B'}C \choose\overline{B'}C}-\overline{A}(-t)F_{1}(\varepsilon;B,B';C;x,y)\notag\\
&~+(q-1)\overline{A}(-t) C(-1)\delta(B'\overline{C})+(q-1)AC(-1)B'\overline{C}(x)\delta(1-y-t)\delta(AB'\overline{C})\\
&~+B(-1)\overline{A}(-x){AB'\overline{C} \choose A\overline{B}}A\overline{C}(y)\overline{A}\overline{B'}C(1-t-y)B'(1-t)\\
&~+BC(-1)\overline{A}(x)A\overline{C}(y)\overline{A}\overline{B'}C(1-y){}_{2}F_1 \left(\begin{matrix}
A, AB'\overline{C} \\
A\overline{B} \end{matrix}
\bigg| \frac{ty}{x(1-y)}\right).
\end{align*}
\end{thm}
\noindent{\it Proof.} It is easily seen from \eqref{p2-2} that
\begin{align}
\label{t41-1}\sum_{\theta}{A\theta\choose \theta}{AB'\overline{C}\theta\choose A\theta}\theta(t/(1-y))&=(q-1){}_{2}F_1 \left(\begin{matrix}
A, AB'\overline{C} \\
 A \end{matrix}
\bigg|t/(1-y)\right)\\
&=(q-1){AB'\overline{C}\choose A}AB'\overline{C}(1-y)\overline{AB'}C(1-y-t)\notag\\
&~-(q-1)A(1-y)\overline{A}(-t)+(q-1)^2 A(-1)\delta(1-y-t)\delta(AB'\overline{C}),\notag\\
  {}_{2}F_1 \left(\begin{matrix}
A\lambda, A\lambda\mu \\
 A\lambda \end{matrix}
\bigg|t\right)&={A\lambda\mu\choose A\lambda}\overline{A}\overline{\lambda}\overline{\mu}(1-t)-\overline{A}\overline{\lambda}(-t). \label{t41-2}
\end{align}
From \eqref{p2-3} and \eqref{t41-2} we know that
\begin{align*}
  &{}_{3}F_2 \left(\begin{matrix}
A, A\lambda, A\lambda\mu \\
A, A\lambda \end{matrix}
\bigg|t\right)\\
&={A\lambda\choose A}{}_{2}F_1 \left(\begin{matrix}
A\lambda, A\lambda\mu \\
 A\lambda \end{matrix}
\bigg|t\right)-\overline{A}(-t){\lambda\mu\choose \lambda}+(q-1)A(-1)\overline{A}\overline{\lambda}(t)\overline{\mu}(1-t)\delta(A\lambda)\\
&={A\lambda\choose \lambda}{A\lambda\mu\choose A\lambda}\overline{A}\overline{\lambda}\overline{\mu}(1-t)-{A\lambda\choose \lambda}\overline{A}\overline{\lambda}(-t)-\overline{A}(-t){\lambda\mu\choose \lambda}\\
&~+(q-1)A(-1)\overline{A}\overline{\lambda}(t)\overline{\mu}(1-t)\delta(A\lambda).
\end{align*}
Thus, by \eqref{f2}--\eqref{f3},
\begin{align}
  &\frac{1}{(q-1)^2}\sum_{\theta,\lambda,\mu}{A\theta\choose \theta}{A\lambda\theta\choose A\theta}{A\lambda\mu\theta\choose A\lambda\theta}{\overline{B}\overline{B'}C \choose \overline{B'}C\lambda}{\overline{B'}C\lambda\choose C\lambda\mu}\lambda(-x)\mu(-y)\theta(t)\label{t41-3}\\
&=\frac{1}{(q-1)^2}\sum_{\lambda,\mu}{\overline{B}\overline{B'}C \choose \overline{B'}C\lambda}{\overline{B'}C\lambda\choose C\lambda\mu}\lambda(-x)\mu(-y)\sum_{\theta}{A\theta\choose \theta}{A\lambda\theta\choose A\theta}{A\lambda\mu\theta\choose A\lambda\theta}\theta(t)\notag\\
&=\frac{1}{q-1}\sum_{\lambda,\mu}{\overline{B}\overline{B'}C \choose \overline{B'}C\lambda}{\overline{B'}C\lambda\choose C\lambda\mu}\lambda(-x)\mu(-y){}_{3}F_2 \left(\begin{matrix}
A, A\lambda, A\lambda\mu \\
A, A\lambda \end{matrix}
\bigg|t\right)\notag\\
&=\frac{\overline{A}(1-t)}{q-1}\sum_{\lambda,\mu}{A\lambda\choose \lambda}{A\lambda\mu\choose A\lambda}{\overline{B}\overline{B'}C \choose \overline{B'}C\lambda}{\overline{B'}C\lambda\choose C\lambda\mu}\lambda\left(\frac{x}{t-1}\right)\mu\left(\frac{y}{t-1}\right)\notag\\
&~-\frac{\overline{A}(-t)}{q-1}\sum_{\lambda,\mu}{A\lambda\choose \lambda}{\overline{B}\overline{B'}C \choose \overline{B'}C\lambda}{\overline{B'}C\lambda\choose C\lambda\mu}\lambda\left(\frac{x}{t}\right)\mu(-y)\notag\\
&~-\frac{\overline{A}(-t)}{q-1}\sum_{\lambda,\mu}{\lambda\mu\choose \lambda}{\overline{B}\overline{B'}C \choose \overline{B'}C\lambda}{\overline{B'}C\lambda\choose C\lambda\mu}\lambda(-x)\mu(-y)\notag\\
&~+B(-1)\overline{A}(-x){AB'\overline{C} \choose A\overline{B}}\sum_{\mu}{B'\mu\choose \overline{A}C\mu}\mu\left(\frac{y}{1-t}\right).\notag
\end{align}
From Theorem \ref{t2-1} we see that
\begin{align}
  &\sum_{\lambda,\mu}{A\lambda\choose \lambda}{A\lambda\mu\choose A\lambda}{\overline{B}\overline{B'}C \choose \overline{B'}C\lambda}{\overline{B'}C\lambda\choose C\lambda\mu}\lambda\left(\frac{x}{t-1}\right)\mu\left(\frac{y}{t-1}\right)\label{t41-4}\\
&=(q-1)^2F_{1}\left(A;B,B';C;\frac{x}{1-t}, \frac{y}{1-t}\right)\notag\\
&~-(q-1)^2C(-1)A(1-t){AB'\overline{C}\choose B'\overline{C}}B'\overline{C}(x)\overline{A}\overline{B'}C(1-t-y)\notag\\
&~-(q-1)^2BC(-1)AB'(1-t)\overline{A}(x)A\overline{C}(y){AB'\overline{C}\choose A\overline{B}}\overline{A}\overline{B'}C(1-t-y).\notag
\end{align}
By \eqref{f1} and \cite [(2.11)]{Gr},
\begin{align}
 &\sum_{\lambda,\mu}{A\lambda\choose \lambda}{\overline{B}\overline{B'}C \choose \overline{B'}C\lambda}{\overline{B'}C\lambda\choose C\lambda\mu}\lambda\left(\frac{x}{t}\right)\mu(-y)\label{t41-5}\\
&=B'(-1)\sum_{\lambda}{A\lambda\choose \lambda}{B\lambda \choose \overline{B'}C\lambda}\lambda\left(\frac{x}{t}\right)\sum_{\mu}{B'\mu\choose C\lambda\mu}\mu(y)\notag\\
&=(q-1)B'(-1)\overline{C}(y)\overline{B'}C(1-y)\sum_{\lambda}{A\lambda\choose \lambda}{B\lambda \choose \overline{B'}C\lambda}\lambda\left(\frac{x(1-y)}{ty}\right)\notag\\
&=(q-1)^2B'(-1)\overline{C}(y)\overline{B'}C(1-y){}_{2}F_1 \left(\begin{matrix}
A, B \\
\overline{B'}C \end{matrix}
\bigg| \frac{x(1-y)}{ty}\right). \notag
\end{align}
It can be deduced from \eqref{f2}, \eqref{f4} and \cite [(2.11)]{Gr} that
\begin{align*}
  \sum_{\mu}{\mu\choose \varepsilon}{\overline{B'}C\choose C\mu}\mu(-y)&=-C(-1)\sum_{\mu}{B'\mu\choose C\mu}\mu(y)+(q-1){\overline{B'}C\choose C}\\
&=-(q-1)\overline{C}(-y)\overline{B'}C(1-y)+(q-1){\overline{B'}C\choose C}.
\end{align*}
This combines \eqref{f4} to give
\begin{align*}
&\sum_{\lambda,\mu}{\lambda\choose \lambda}{\lambda\mu\choose \lambda}{\overline{B}\overline{B'}C \choose\overline{B'}C\lambda}{\overline{B'}C\lambda\choose C\lambda\mu}\lambda(-x)\mu(-y)\\
&=(q-1){\overline{B}\overline{B'}C \choose\overline{B'}C}\sum_{\mu}{\mu\choose \varepsilon}{\overline{B'}C\choose C\mu}\mu(-y)-\sum_{\lambda,\mu}{\lambda\mu\choose \lambda}{\overline{B}\overline{B'}C \choose\overline{B'}C\lambda}{\overline{B'}C\lambda\choose C\lambda\mu}\lambda(-x)\mu(-y)\\
&=(q-1)^2 {\overline{B'}C\choose C}{\overline{B}\overline{B'}C \choose\overline{B'}C}-(q-1)^2 \overline{C}(-y)\overline{B'}C(1-y){\overline{B}\overline{B'}C \choose\overline{B'}C}\\
&~-\sum_{\lambda,\mu}{\lambda\mu\choose \lambda}{\overline{B}\overline{B'}C \choose\overline{B'}C\lambda}{\overline{B'}C\lambda\choose C\lambda\mu}\lambda(-x)\mu(-y).
\end{align*}
From Theorem \ref{t2-1} we have
\begin{align*}
 &\sum_{\lambda,\mu}{\lambda\choose \lambda}{\lambda\mu\choose \lambda}{\overline{B}\overline{B'}C \choose \overline{B'}C\lambda}{\overline{B'}C\lambda\choose C\lambda\mu}\lambda(-x)\mu(-y)\\
&=(q-1)^2F_{1}(\varepsilon;B,B';C;x,y)+(q-1)^2C(-1)B'\overline{C}\left(\frac{x}{1-y}\right)\\
&~-(q-1)^3 C(-1)\delta(B'\overline{C})-(q-1)^2BC(-1)\overline{C}(y)\overline{B'}C(1-y){B'\overline{C}\choose \overline{B}}.\\
.
\end{align*}
So we deduce from the above two identities that
\begin{align}
&\sum_{\lambda,\mu}{\lambda\mu\choose \lambda}{\overline{B}\overline{B'}C \choose\overline{B'}C\lambda}{\overline{B'}C\lambda\choose C\lambda\mu}\lambda(-x)\mu(-y)\label{t41-6}\\
&=(q-1)^2 {\overline{B'}C\choose C}{\overline{B}\overline{B'}C \choose\overline{B'}C}-(q-1)^2F_{1}(\varepsilon;B,B';C;x,y)\notag\\
&~-(q-1)^2C(-1)B'\overline{C}\left(\frac{x}{1-y}\right)+(q-1)^3 C(-1)\delta(B'\overline{C}).\notag
\end{align}
It follows from \cite [(2.11)]{Gr} that
\begin{equation}\label{t41-7}
\sum_{\mu}{B'\mu\choose \overline{A}C\mu}\mu\left(\frac{y}{1-t}\right)=(q-1)A\overline{C}(y)\overline{A}\overline{B'}C(1-t-y)B'(1-t).
\end{equation}
By Theorem \ref{t2-1} and \eqref{f2}--\eqref{f3},
\begin{align*}
  &\sum_{\theta}{A\theta\choose\theta}F_{1}(A\theta;B,B';C;x,y)\theta(t)\\
&=\frac{1}{(q-1)^2}\sum_{\theta,\lambda,\mu}{A\theta\choose \theta}{A\lambda\theta\choose A\theta}{A\lambda\mu\theta\choose A\lambda\theta}{\overline{B}\overline{B'}C \choose \overline{B'}C\lambda}{\overline{B'}C\lambda\choose C\lambda\mu}\lambda(-x)\mu(-y)\theta(t)\\
&~+C(-1)B'\overline{C}(x)\overline{A}\overline{B'}C(1-y)\sum_{\theta}{A\theta\choose \theta}{AB'\overline{C}\theta\choose A\theta}\theta(t/(1-y))\\
&~+BC(-1)\overline{A}(x)A\overline{C}(y)\overline{A}\overline{B'}C(1-y)\sum_{\theta}{A\theta\choose \theta}{AB'\overline{C}\theta\choose A\overline{B}\theta}\theta\left(\frac{yt}{x(1-y)}\right).
\end{align*}
Substituting \eqref{t41-4}--\eqref{t41-7} into \eqref{t41-3}, then applying \eqref{t41-1} and \eqref{t41-3} in the above identity and cancelling some terms, we can easily obtain the result. This ends the proof of Theorem \ref{t4-1}.\qed
\section*{Acknowledgement}
This  work was  supported by the Natural Science Basic Research Plan in Shaanxi Province of China (No. 2017JQ1001), the Initial Foundation for Scientific Research  of Northwest A\&F University (No. 2452015321) and the Fundamental Research Fund of Northwest A\&F University (No. 2452017170).  

\end{document}